\documentclass[11pt]{article}
\usepackage{amssymb,amsmath,mathrsfs,txfonts,graphicx,color}
\graphicspath{{figs/}}
\usepackage[numbers,sort&compress]{natbib}
\allowdisplaybreaks
\def\pd#1#2{\frac{\partial#1}{\partial#2}}
\let\oldsection\section
\renewcommand\section{\setcounter{equation}{0}\oldsection}

\newtheorem{theorem}{\indent Theorem}[section]
\newtheorem{lemma}{\indent Lemma}[section]

\newtheorem{corollary}{\indent Corollary}[section]

\begin{document}

\title{\LARGE
Propagation speed of degenerate diffusion equations \\
with time delay}
\author{
Tianyuan Xu$^{a}$, Shanming Ji$^{a,}$\thanks{Corresponding author, email:jism@scut.edu.cn}, \
Ming Mei$^{b,c}$, Jingxue Yin$^d$
\\
\\
{ \small \it $^a$School of Mathematics, South China University of Technology}
\\
{ \small \it Guangzhou, Guangdong, 510641, P.~R.~China}
\\
{ \small \it $^b$Department of Mathematics, Champlain College Saint-Lambert}
\\
{ \small \it Quebec,  J4P 3P2, Canada, and}
\\
{ \small \it $^c$Department of Mathematics and Statistics, McGill University}
\\
{ \small \it Montreal, Quebec,   H3A 2K6, Canada}
\\
{ \small \it $^d$School of Mathematical Sciences, South China Normal University}
\\
{ \small \it Guangzhou, Guangdong, 510631, P.~R.~China}
}
\date{}

\maketitle

\begin{abstract}
We are concerned with a class of degenerate diffusion equations with time delay describing population dynamics with age structure.
In our recent study [{\em Nonlinearity}, 33 (2020), 4013--4029], we
established the existence and uniqueness of critical traveling wave for the time-delayed degenerate diffusion equations, and obtained
the reducing mechanism of  time delay on critical wave speed. In this paper,
we now are able to show the asymptotic spreading speed and its coincidence with the critical wave speed $c^*(m,r)$ of sharp wave,
and prove that the initial perturbation or the boundary of the compact support of the solution propagates at the critical wave speed $c^*(m,r)$ for the time-delayed degenerate diffusion equations.
Remarkably, different from the existing studies related to spreading speeds, the time delay and the degenerate diffusion lead to some essential difficulties in the analysis of the spreading
speed, because the time-delay makes the critical speed of traveling waves slow down, and the degenerate diffusion causes the loss of regularity for the solutions.
By a phase transform technique combined with the monotone method, we can determine the asymptotic spreading speed.
Furthermore, we propose a brand-new sharp-profile-based difference scheme to handle large variation
of degenerate diffusion $(u^m)_{xx}$ near the sharp edge and
carry out some numerical simulations which perfectly confirm our theoretical results.
\end{abstract}


\section{Introduction}

We consider the following degenerate diffusion equation with time delay
\begin{align}\label{eq-main}
\begin{cases}
\displaystyle
\pd u t=\Delta u^m-d(u)+b(u(t-r,x)),\quad &x\in \mathbb R^n,~t>0,
\\[1mm]
u(s,x)=u_0(s,x), \quad &x\in\mathbb R^n,~t\in[-r,0],
\end{cases}
\end{align}
where the spatial dimension $n\ge1$,
$u(t,x)$ denotes the total mature population of the species at location $x$ and time $t>0$,
$r\ge0$ is the maturation time, $b(u(t-r,x))$ is the birth function,
$d(u)$ is the death rate function.
The equation \eqref{eq-main}
describes the population dynamics of single species with age-structure and density-dependent
diffusion ($m>1$).

The positive density dependence on emigration has been predicted by many theoretical models and confirmed empirically in various species with age structure \cite{Denno,Murry,Nowicki}.
Time delays arise from the passage through sequential demographic phases in the species' life cycle \cite{Sulsky}.
Hence, time delay and nonlinear dispersal are both inextricably mechanism in population dynamics.
For example, butterflies, the most popular age-structured species for dispersal studies,  have
been reported increased emigration at high population densities \cite{Nowicki}.
Positive density-dependent dispersal is particularly
beneficial for females, because it allows laying the egg-load in less crowded
patches to help their offspring avoiding severe intra-specific
competition in the larval period.

In the linear diffusion case without time delay, (i.e., $r=0$ and $m=1$),
the equation \eqref{eq-main} is reduced to the Fisher and Kolmogorov-Petrovsky-Piskunov (KPP) equation \cite{Fisher37,Kolmogorov}.
It is well known that there exists a critical (minimal admissible) wave speed $c^*=2\sqrt{b'(0)-d'(0)}>0$
(under certain conditions on the functions $b(\cdot)$ and $d(\cdot)$) for all the traveling waves
connecting the two constant equilibria $0$ and $\kappa>0$,
and the level set $\Gamma_\varepsilon(t):=\{x\in\mathbb R^n;u(t,x)=\varepsilon\}$ with $\varepsilon\in(0,\kappa)$
asymptotically propagates at the same speed $c^*$  \cite{Aronsonw78,Weinberger82}.
It was Thieme and Zhao \cite{Thieme} who first established the theory of asymptotic spreading speed for a
large class of  nonlinear integral equations, which
covers many time-delayed reaction and diffusion equations with linear diffusion (i.e., $r>0$ and $m=1$).
Liang and Zhao \cite{LiangX07} further
developed the theory of spreading speeds to both discrete and continuous time
monotone semiflows and investigated the application to a time-delayed evolution equation.
In a series of works (e.g. \cite{LiZhao,ZhaoCV20,ZhangZhaoJDE20}), Zhao and his collaborators
investigated the traveling waves and spreading speeds of population dynamics model with nonlocal dispersal.
Studies of the coincidence of the spreading speed with the critical wave speed for
various evolution systems with linear diffusion systems can also be found in
\cite{BerestyckiJAMS,Chern-Mei-Zhang-Yang-JDE15,LiB05,LiangX06,LiangX07,LiangX10,Mei-LinJDE09,Thieme,zhangzhao}.

When the degenerate diffusion is included, the  system can be used to describe biological population dynamics
with density-dependent dispersal; see for instance \cite{Gurney75,Murry}.
An interesting peculiarity of degenerate diffusion is the appearance of sharp type waves at the asymptotic speed \cite{BenguriaPRE,PME,VazquezSmooth06}.
For the case without time delay ($r=0$), traveling wave solutions have
been found by several authors \cite{Aronson80Density,Aronson86PME,Pablo,Gilding,Garduno}.
Medvedev et al. \cite{EJAM} proved that the slowest traveling wave in the family yields the asymptotic speed
of the propagation of disturbances in a class of degenerate Fisher-KPP equations.
In recent works \cite{AudritoDCDS2019,AudritoJDE2017,AudritoJDE2020},
more general cases of doubly nonlinear diffusion are considered,
which includes both porous medium and $p$-Laplacian models.

An increasing attention has been paid to degenerate diffusion equations with time delay
in order to study the effects of degenerate diffusion and time delay on
the evolutionary behavior of biological systems with age structure (see \cite{HJMY,JDE18,Non20,JDE20}).
The existence of smooth monotone fronts
for equations \eqref{eq-main} with small time delay
was proved by Huang et al.~\cite{HJMY} based on a perturbation approach.
%Later then,
%we \cite{JDE18} proved the existence of non-critical smooth monotone traveling wave solutions for any time delay by the upper and lower solutions method combining with the viscosity vanishing method.
%In our recent work \cite{Non20}, we established the existence of
%critical sharp type traveling wave and investigate how the time delay affects the
%propagation mechanism of fronts.
In our recent works \cite{JDE18,Non20}, we proved that the time-delayed degenerate diffusion equation
\eqref{eq-main} admits a unique sharp type (semi-compactly supported) traveling wave
$\phi(x+c^*t)$ for the one dimensional case,
which corresponds to the minimal admissible (critical) wave speed $c^*=c^*(m,r)$.
Moreover, the time delay slows down the minimal wave speed,
i.e., $c^*(m,r)<c^*(m,0)$ for $r>0$.
However, the asymptotic speeds of
spreading for solutions with compactly supported initial data and the coincidence with the critical
wave speed of sharp traveling wave still remain open.
In this paper, we shall answer those unsettled questions on the spreading properties.
%show that they coincide with the asymptotic speeds of
%spread for solutions with initial functions having compact supports.

Time delay and degenerate diffusion
lead to essential difficulties in the analysis of the spreading speed of \eqref{eq-main}.
In the absence of time delay, the maximum principle and phase plane analysis proposed by
Aronson and Weinberger \cite{AronsonW75,Aronsonw78} yield conclusions about the asymptotic propagation speed
of the linear diffusion and similar equations.
This method was extended to cover degenerate diffusion equations of general Fisher-KPP sources without time delay
by Medvedev et al.~in \cite{EJAM},
where all the trajectories in the phase plane are determined and correspond to special
upper and lower solutions.
However, time delay changes the situation dramatically.
It is shown in \cite{Non20} that the time delay reduces the critical wave speed $c^*(m,r)$
and the speed is not characterized by the classical phase plane analysis method.
In order to construct upper and lower solutions with compact (or semi-compact) supports
and with propagating speed approaching $c^*(m,r)$,
we employ a new phase transform technique developed in \cite{Non20}
and utilize the monotone dependence in the phase space,
see the phase comparison principle Lemma \ref{le-monotone} in this paper.
Especially, we  need to treat the lower solutions technically in two ways:
first, we show that the support of $u(t,\cdot)$ expands
to including any given compact subset for large time;
and secondly, the value of
$u(t,x)$ within given compact subset grows up as time increases.

The purpose of this paper is to study the propagation speed for the equation \eqref{eq-main}.
%Regarding the model \eqref{eq-main} with both effects of degeneracy of diffusion and time-delay, based on the concerns of populations with age-structure, a brief derivation of the model will be made in Section 2.
The main issue of the paper is to show that the initial perturbation or the boundary of the compact support of the solution
propagates at the critical wave speed $c^*(m,r)$ for the time-delayed degenerate diffusion equation \eqref{eq-main}.
Different from the existing studies related to spreading speeds, the time delay and the degenerate diffusion lead to some essential difficulties in the analysis of the spreading
speed, because the time-delay makes the critical speed of traveling waves slow down, and the degenerate diffusion causes the loss of regularity for the solutions. This main result will be proved in Section 2.
Section 3 is devoted to numerical computations. Since the variation of the degenerate diffusion $(u^m)_{xx}$ may be large near
the sharp edge, the traditional numerical schemes are failed in the case with moving sharp-edge. Here, we propose a brand new numerical algorithm, called the sharp-profile-based difference scheme,
and use this scheme to carry out some numerical simulations in different cases, which are accurate and stable in the sense of numerical performance,  and also perfectly demonstrate our theoretical results.

In the end of this section, we are going to state our main theorems on the propagation speed for the time-delayed degenerate-diffusion equations \eqref{eq-main}.
A function $u(t,x)$ is said to be compactly supported for $t\in[t_1,t_2]$,
if $\mathrm{supp}\,u(t,\cdot)$ is compact for $t\in[t_1,t_2]$.
For the sake of convenience, we define the half space divided by the hyperplane through a point $x_0$
that has normal vector $\nu$ as
\begin{equation} \label{eq-semi}
\Pi(x_0,\nu)
:=\{x\in\mathbb R^n;(x-x_0)\cdot\nu\ge0\}.
\end{equation}
A function $u(t,x)$ is said to be semi-compactly supported for $t\in[t_1,t_2]$,
if
\begin{equation*} \label{eq-partial-supp}
\mathrm{supp}\,u(t,\cdot)\subset \Pi(x(t),\nu(t))
=\{x\in\mathbb R^n;(x-x(t))\cdot\nu(t)\ge0\}
\end{equation*}
for $t\in[t_1,t_2]$ with some points $x(t)$ and vectors $\nu(t)$.
For any variable $s$, the positive value of $s$, denoted by $s_+$,
is defined as $s_+=\max\{s,0\}$.

Throughout this paper, we assume that the functions $d(s)$ and $b(s)$ satisfy the following
conditions:
\begin{align} \nonumber
&\text{There exist~}u_1=0, u_2>0 \text{~such that~}
d, b\in C^2([0,+\infty)),
d(0)=b(0)=0,d(u_2)=b(u_2),
\\ \label{eq-H}
&(b(s)-d(s))\cdot(s-u_2)<0, \forall s\in(0,+\infty)\backslash\{u_2\},
b'(0)>d'(0)\ge 0, d'(s)\ge 0, b'(s)\ge0.
\end{align}
Here, $u_1=0$ and $\kappa:=u_2>0$ are two constant equilibria of \eqref{eq-main},
and the functions $b(u)$, $d(u)$ are both non-decreasing.

We first recall the properties of the sharp type traveling wave obtained in our previous study \cite{Non20}.

\begin{theorem}[\cite{Non20}] \label{le-sharp}
For any $m>1$ and $r\ge0$, the time-delayed degenerate diffusion equation \eqref{eq-main}
admits a unique (up to shift) sharp type (semi-compactly supported) traveling wave
$u^*(t,x)=\phi^*(x\cdot\nu+c^*t)$ with a unique speed $c^*=c^*(m,r)>0$ and unit vector $\nu$ such that:

$\mathrm(i)$ $\mathrm{supp}\,u^*(t,\cdot)=\Pi(-c^*t\nu,\nu)$, $\phi^*(\xi)$ is monotone increasing
and $\phi^*(+\infty)=\kappa$;

$\mathrm(ii)$ $c^*(m,r)$ is the minimal admissible (or critical) traveling wave speed;

$\mathrm(iii)$ $c^*(m,r)<c^*(m,0)$, i.e., the time delay slows down the critical traveling wave speed.
\end{theorem}

If the initial data $u_0(s,x)$ are compactly supported (or semi-compactly supported) in a direction $\nu\in\mathbb S^{n-1}$,
say $\mathrm{supp}\,u_0(s,\cdot)\subset \Pi(x_0,\nu)$ for example,
it is expected that the solution $u(t,x)$ is compactly supported (or semi-compactly supported)
in this direction for all the time
since the diffusion of the equation \eqref{eq-main} is degenerate.
The solution $u(t,x)$ together with its support $\mathrm{supp}\,u(t,\cdot)$ expands toward opposite $\nu$ direction.
In what follows, we always assume that $c^*(m,r)>0$ is the critical (minimal admissible) wave speed of traveling waves
and also is the wave speed of the unique sharp type (semi-compactly supported)
traveling wave of \eqref{eq-main} shown by Theorem \ref{le-sharp}.
Now we present the following large time propagation speed of the solution
with compact (or semi-compact) support.

\begin{theorem} \label{th-evolution}
Let $u(t,x)$ be the solution of \eqref{eq-main} with initial data
$u_0(s,x)$ for $s\in[-r,0]$ satisfying
\begin{equation} \label{eq-initial-semi}
\mathrm{supp}\,u_0(s,\cdot)\subset \Pi(x_0,\nu),
\qquad u_0(s,x)\ge \phi_0((x-x_0)\cdot\nu),
\end{equation}
where $\phi_0(\eta)$ with $\eta=(x-x_0)\cdot\nu$
is a non-negative, continuous and non-trivial function.
Then for any $0<c_1<c^*(m,r)<c_2$, and any $0<\kappa_1<\kappa<\kappa_2$,
there exist a time $T=T(c_1,c_2,\kappa_1,\kappa_2)>0$
and two functions $\phi_1(\cdot)$ and $\phi_2(\cdot)$ with semi-compact supports, such that
\begin{equation} \label{eq-propagating}
\phi_1((x-x_0)\cdot\nu+c_1t)\le u(t,x)\le \phi_2((x-x_0)\cdot\nu+c_2t), \quad \forall (x-x_0)\cdot\nu\le0, ~t\ge T,
\end{equation}
and
\begin{align*}
\mathrm{supp}\,\phi_1(\cdot)=\mathrm{supp}\,\phi_2(\cdot)= \Pi(x_0,\nu), \quad
\lim_{\eta\to+\infty}\phi_1(\eta)\ge\kappa_1, \quad \lim_{\eta\to+\infty}\phi_2(\eta)\le\kappa_2.
\end{align*}
Therefore,
\begin{equation} \label{eq-support}
\Pi(x_0-c_1t,\nu)\cap\Pi(x_0,-\nu)\subset \mathrm{supp}\,u(t,\cdot)\subset \Pi(x_0-c_2t,\nu), \quad t\ge T.
\end{equation}
\end{theorem}

\begin{corollary}
Under the condition of Theorem \ref{th-evolution}, there holds
$$
\lim_{t\to+\infty} \, \inf_{x\in\mathbb R^n, \ u(t,x)>0}\frac{(x-x_0)\cdot\nu}{t}=-c^*(m,r),
$$
and
$$
\lim_{t\to+\infty} \, \sup_{x\in \Pi(x_0,-\nu), \ u(t,x)=0}\frac{(x-x_0)\cdot\nu}{t}=-c^*(m,r).
$$
Additionally, for any $x\in\mathbb R^n$, it holds
\begin{equation} \label{eq-large-speed}
\lim_{t\to+\infty}u(t,x-ct\nu)=
\begin{cases}
\kappa, \quad &\text{if~} c<c^*(m,r), \\
0, \quad &\text{if~} c>c^*(m,r).
\end{cases}
\end{equation}
\end{corollary}

Theorem \ref{th-evolution} implies that, the solution with non-trivial and non-negative
initial data that are compactly supported in one direction
(or compactly supported for the one dimensional case)
propagates at the same speed as the sharp type traveling wave,
which is the minimal wave speed of all the traveling waves.
Particularly for the one dimensional case, any compactly supported initial perturbation remains compact
and the boundary propagates at speed $c^*(m,r)$.

%The paper is organized as follows.
%In Section 2, we prove the local propagation behavior,
%including the persistence, immediate positivity, waiting time, and propagation speed.
%The large time speed of propagation is shown in Section 3.

\section{Proof of the main result}

We are interested in the propagation speed of the solutions.
We assume that $u_0(s,x)$ for $s\in[-r,0]$ is non-trivial, non-negative,
bounded and continuous,
therefore the Cauchy problem \eqref{eq-main} can be solved step by step
such that $u(t,x)$ is non-negative, bounded and continuous on $x\in\mathbb R^n$
and $t\in[-r,+\infty)$ (see \cite{PME} for example).
Moreover, the comparison principle holds for the Cauchy problem \eqref{eq-main}
and the initial boundary value problem on bounded domain
since the time-delayed source $b(u(t-r,x))$ is monotone increasing
with respect to $u(t-r,x)$.
We shall prove that for large time scale the average speed of propagation
is consistent with the sharp type traveling wave speed.

The sharp type traveling wave is the typical solution that is semi-compactly supported
and propagates at a positive and finite speed $c^*(m,r)$.
In order to show the large time propagation speed of general
solutions with compact (or semi-compact) supports, we need to
construct upper and lower solutions with compact (or semi-compact) supports
and propagates at speed near $c^*(m,r)$.
The case of $r=0$ (no time delay) is proved by the phase-plane analysis,
where all the trajectories are determined and correspond to special
upper and lower solutions, see \cite{EJAM} and the references therein.
Here for the time-delayed case ($r>0$), we need to employ a new phase transform approach
developed in \cite{Non20}.

Consider the following ``traveling wave'' type special function
defined for any $c>0$ and $k\ge\kappa$
\begin{equation} \label{eq-tw}
u_c^k(t,x):=\phi_c^k(x\cdot\nu+ct)
=\phi_c^k(\xi), \quad \text{with~} \xi:=x\cdot\nu+ct,
\end{equation}
such that
\begin{equation} \label{eq-tw2}
\phi_c^k(\xi)=0, \quad\forall\xi\le0,
\qquad \phi_c^k(\xi)\in(0,k),
\quad \forall \xi\in(0,\xi_c^k),
\end{equation}
for some $\xi_c^k\in(0,+\infty]$.
Note that the sharp type traveling wave $\phi^*(\xi)$ in Theorem \ref{le-sharp}
is a traveling wave type special function $\phi_{c^*}^{\kappa}(\xi)$
corresponding to the critical traveling wave speed $c^*=c^*(m,r)$.
As proved in \cite{Non20}, $c^*$ is the unique speed
and $\phi^*$ is the unique function
such that $\phi^*$ satisfies \eqref{eq-tw}, \eqref{eq-tw2},
and the time-delayed degenerate diffusion equation \eqref{eq-main}$_1$.
Therefore, we only expect that $\phi_c^k$ is a local solution
of \eqref{eq-main}$_1$ near the boundary $\xi=0$ for $c\ne c^*$.
Actually, this is solved through a delayed iteration scheme as follows.

The traveling wave type function $\phi_c^k(\xi)$
defined by \eqref{eq-tw} and \eqref{eq-tw2}
is a local solution of \eqref{eq-main}$_1$ on $\xi\in(-\infty,\xi_c^k)$
with $\xi_c^k\in(0,+\infty]$ if
\begin{equation} \label{eq-tw-local}
\begin{cases}
c(\phi_c^k)'(\xi)=((\phi_c^k)^m(\xi))''
-d(\phi_c^k(\xi))+b(\phi_c^k(\xi-cr)),
\quad \xi\in(-\infty,\xi_c^k), \\
\phi_c^k(\xi)=0, \quad \forall\xi\le0, \qquad \phi_c^k(\xi)\in(0,k),
\quad \forall \xi\in(0,\xi_c^k).
\end{cases}
\end{equation}

The local solvability of the degenerate second order differential equation
\eqref{eq-tw-local} is proved in the following lemma.

\begin{lemma} \label{le-solve}
For any $c>0$ and $k\ge \kappa$, the degenerate problem \eqref{eq-tw-local}
admits a unique local solution $\phi_c^k(\xi)$ on $(-\infty,\xi_c^k)$
with $\xi_c^k\in(0,+\infty]$
(we may assume that $(-\infty,\xi_c^k)$ is the maximal existence interval)
such that
$$
\phi_c^k(\xi)=\Big(\frac{(m-1)c}{m}\xi\Big)_+^\frac{1}{m-1}+o(|\xi|^\frac{1}{m-1}),
\quad \xi\to0.
$$
Moreover, (a) $\phi_c^k(\xi)$ is strictly increasing on $(0,\xi_c^k)$
and $\phi_c^k(\xi_c^k)=k$,
or (b) $\phi_c^k(\xi)$ is not strictly increasing on $(0,\xi_c^k)$ and $\phi_c^k(\xi_c^k)=0$.
For the case (b), there holds $\xi_c^k<+\infty$ and there exists a $\hat\xi_c^k\in(0,\xi_c^k)$
such that $\phi_c^k(\xi)$ is strictly increasing on $(0,\hat\xi_c^k)$
and decreasing on $(\hat\xi_c^k,\xi_c^k)$.
\end{lemma}
{\it\bfseries Proof.}
Note that $\phi_c^k$ is semi-compactly supported and the time-delayed source
function $b(\phi_c^k(\xi-cr))=0$ if $\xi\le cr$.
Therefore, \eqref{eq-tw-local} is locally reduced to the following equation
\begin{equation} \label{eq-tw-local-cr}
\begin{cases}
c\phi_c'(\xi)=(\phi_c^m(\xi))''-d(\phi_c(\xi)),
\quad \xi\in(-\infty,cr), \\
\phi_c(\xi)=0, \quad \forall\xi\le0, \qquad \phi_c(\xi)>0,
\quad \forall \xi\in(0,cr),
\end{cases}
\end{equation}
whose unique solvability is proved in \cite{Non20,JDE20}
and the solution is denoted by $\phi_c$.
Lemma 3.1 (and its proof) in \cite{JDE20} shows
the locally asymptotical behavior of $\phi_c(\xi)$
near zero and the strictly increasing monotonicity of $\phi_c(\xi)$ in $(0,cr)$.
If $\phi_c(cr)\ge k$, then $\xi_c^k=\sup\{\xi\in(0,cr);\phi_c(\xi)<k\}$.
If $\phi_c(cr)< k$, then we solve \eqref{eq-tw-local} on $(cr,2cr)$ as
\begin{equation} \label{eq-tw-local-2cr}
\begin{cases}
c\phi_c'(\xi)=(\phi_c^m(\xi))''
-d(\phi_c(\xi))+b(\phi_c(\xi-cr)),
\quad \xi\in(cr,2cr), \\
\phi_c(cr), \phi_c'(cr) \text{~are~determined~from~left~side}.
\end{cases}
\end{equation}
The problem \eqref{eq-tw-local-2cr} is locally solved near $cr$
since $\phi_c(cr)>0$ and $b(\phi_c(\xi-cr))$ is already known from \eqref{eq-tw-local-cr}.
Then three cases may happen:

(i) $\phi_c(\xi)$ is not strictly increasing on whole $(cr,2cr)$,
which means there exists a $\xi_0\in(cr,2cr)$ such that $\phi_c'(\xi_0)\le0$.
We employ Lemma 3.5 in \cite{Non20} to derive that
$\phi_c(\xi)$ is always decreasing after $\xi_0$ until reaching zero for $\xi>\xi_0$.
If $\phi_c(\xi)>0$ for all $\xi\in(cr,2cr)$, then we solve \eqref{eq-tw-local}
further on $(2cr,3cr)$ and the intervals after this
in a similar way as \eqref{eq-tw-local-2cr}
until $\phi_c(\xi_1)=0$ for some
$\xi_1>\xi_0>cr$.
In this case, $\xi_c^k=\sup\{\xi>cr;\phi_c(\xi)>0\}$ and $\phi_c(\xi_c^k)=0$.
The assertion $\xi_c^k<+\infty$ is proved in a similar way as the proof of
Lemma 3.5 in \cite{JDE20}.

(ii) $\phi_c(\xi)$ is strictly increasing on $(cr,2cr)$
and $\phi_c(2cr)\ge k$, then $\xi_c^k=\sup\{\xi\in(cr,2cr);\phi_c(\xi)<k\}$
and $\phi_c(\xi_c^k)=k$.

(iii) $\phi_c(\xi)$ is strictly increasing on $(cr,2cr)$
and $\phi_c(2cr)< k$, then we solve \eqref{eq-tw-local} further on $(2cr,3cr)$
and the intervals after this until (i) or (ii) happens.
Otherwise, $\phi_c(\xi)$ is strictly increasing and \eqref{eq-tw-local} is solved on
$(-\infty,+\infty)$ such that $\xi_c^k=+\infty$ and $\phi_c(\xi_c^k)=k$.
This happens for $c=c^*$ and $k=\kappa$ since $\phi^*(\xi)=\phi_{c^*}^\kappa(\xi)$
is the unique sharp type traveling wave.
$\hfill\Box$

In order to show more precise behavior of $\phi_c^k$,
we employ the following phase transform approach
and formulate phase comparison principle.
For any $c>0$ and $k\ge \kappa$, let $\phi_c^k(\xi)$
be the unique solution of the degenerate problem \eqref{eq-tw-local}
on its maximal existence interval $(-\infty,\xi_c^k)$ with $\xi_c^k\in(0,+\infty]$
as shown in Lemma \ref{le-solve},
and let $\hat\xi_c^k\in(0,\xi_c^k]$ be the largest number
(or equivalently, $(0,\hat\xi_c^k)$ be the maximal interval)
such that $\phi_c^k(\xi)$ is strictly increasing on $(0,\hat\xi_c^k)$.
Define
(here we do not explicitly write down the dependence of $\psi_c(\xi)$ on $k$ for simplicity)
\begin{equation} \label{eq-psi}
\psi_c(\xi):=((\phi_c^k)^m(\xi))'=m(\phi_c^k)^{m-1}(\xi)\cdot(\phi_c^k)'(\xi),
\quad \xi\in(0,\hat\xi_c^k).
\end{equation}
Now we have two functions $\phi_c^k(\xi)$ and $\psi_c(\xi)$ defined for
$\xi\in(0,\hat\xi_c^k)$,
and $\phi_c^k(\xi)$ is strictly increasing on $(0,\hat\xi_c^k)$,
then we can interpret $\psi=\psi_c(\xi)$ as a function of $\phi=\phi_c^k(\xi)$
through the intermediate variable
$$
\xi=(\phi_c^k)^{-1}(\phi), \quad \phi\in(0,\phi_c^k(\hat\xi_c^k)).
$$
That is, we rewrite
\begin{equation} \label{eq-psi-phi}
\tilde\psi_c(\phi):=\psi_c(\xi)=\psi_c((\phi_c^k)^{-1}(\phi)),
\quad \phi\in(0,\phi_c^k(\hat\xi_c^k)).
\end{equation}
A key transform in dealing with the time delay
in the degenerate diffusion equation \eqref{eq-tw-local}
is to rewrite $\phi_c^k(\xi-cr)$ as a function of $\phi=\phi_c^k(\xi)$
depending on $\tilde\psi_c(\phi)$ in a functional way:
\begin{equation} \label{eq-phicr}
\phi_{c,cr}(\phi):=\phi_c^k(\xi-cr)
=\phi_c^k((\phi_c^k)^{-1}(\phi)-cr)
=\inf_{\theta\ge0}\Big\{\int_\theta^\phi\frac{ms^{m-1}}{\tilde\psi_c(s)}\mathrm{d}s\le cr\Big\},
\quad \phi\in(0,\phi_c^k(\hat\xi_c^k)).
\end{equation}

\begin{lemma}[Phase transform] \label{le-phicr}
The functional interpretation \eqref{eq-phicr} is well-defined for the
sharp type functions $\phi_c^k(\xi)$
for $\phi\in(0,\phi_c^k(\hat\xi_c^k))$.
\end{lemma}
{\it\bfseries Proof.}
We divide the proof into two cases.

Case I.
If for some $\phi=\phi_c^k(\xi)$ with $\xi\in(0,\hat\xi_c^k)$
and $\phi\in(0,\phi_c^k(\hat\xi_c^k))$ there holds
$$
\int_0^\phi\frac{ms^{m-1}}{\tilde\psi_c(s)}\mathrm{d}s>cr,
$$
then we rewrite the above integral through the method of substitution
of $s=\phi_c^k(t)$ for $s\in(0,\phi)$ and $t\in(0,\xi)$ to find
$$
cr<\int_0^\phi\frac{ms^{m-1}}{\tilde\psi_c(s)}\mathrm{d}s
=\int_0^\xi\frac{m(\phi_c^k)^{m-1}(t)}{m(\phi_c^k)^{m-1}(t)\cdot(\phi_c^k)'(t)}
(\phi_c^k)'(t)\,\mathrm{d}t=\xi.
$$
Therefore, $\xi-cr>0$ and $\phi_{c,cr}(\phi)=\phi_c^k(\xi-cr)$ is the unique value such that
$$
\int_{\phi_{c,cr}(\phi)}^\phi\frac{ms^{m-1}}{\tilde\psi_c(s)}\mathrm{d}s=cr.
$$

Case II.
If for some $\phi=\phi_c^k(\xi)$ with $\xi\in(0,\hat\xi_c^k)$
and $\phi\in(0,\phi_c^k(\hat\xi_c^k))$ there holds
$$
\int_0^\phi\frac{ms^{m-1}}{\tilde\psi_c(s)}\mathrm{d}s\le cr,
$$
then $\xi-cr\le0$ and $\phi_c^k(\xi-cr)=0$
since $\phi_c^k$ is sharp type such that $\phi_c^k(t)\equiv0$ for all $t\le0$.
$\hfill\Box$

\begin{lemma}[Monotone dependence] \label{le-monotone}
For any $c>0$ and $k\ge \kappa$, let $\phi_c^k(\xi)$
be the unique solution of the degenerate problem \eqref{eq-tw-local}
on its maximal existence interval $(-\infty,\xi_c^k)$ with $\xi_c^k\in(0,+\infty]$
and let $\tilde\psi_c(\phi)$ and $\phi_{c,cr}(\phi)$ be the
phase transform functions defined by \eqref{eq-psi-phi} and \eqref{eq-phicr}.
Then for $c_1>c_2>0$, there holds
\begin{equation} \label{eq-phase-mono}
\tilde\psi_{c_1}(\phi)>\tilde\psi_{c_2}(\phi),
\quad \forall\phi\in(0,\min\{\phi_{c_1}^k(\hat\xi_{c_1}^k),\phi_{c_2}^k(\hat\xi_{c_2}^k\})),
\end{equation}
and
\begin{equation} \label{eq-monotone}
\phi_{c_1}^k(\xi)>\phi_{c_2}^k(\xi), \quad
\forall\xi\in(0,\min\{\hat\xi_{c_1}^k,\hat\xi_{c_2}^k\}).
\end{equation}
\end{lemma}
{\it\bfseries Proof.}
The monotone dependence of $\phi_c^k(\xi)$ with respect to $c$ is
proved in Lemma 3.6 in \cite{Non20}.
The monotone dependence $\tilde\psi_{c_1}(\phi)>\tilde\psi_{c_2}(\phi)$
means that
$$
(\phi_{c_1}^k)'(\xi_1)>(\phi_{c_2}^k)'(\xi_2),
\quad \text{at~where~}\phi_{c_1}^k(\xi_1)=\phi=\phi_{c_2}^k(\xi_2),
$$
or equivalently,
\begin{equation} \label{eq-phase-comp}
(\phi_{c_1}^k)'((\phi_{c_1}^k)^{-1}(\phi))>(\phi_{c_2}^k)'((\phi_{c_2}^k)^{-1}(\phi)).
\end{equation}
In contrast to the comparison between two functions
$\phi_{c_1}^k(\xi)$ and $\phi_{c_2}^k(\xi)$
at the same point $\xi$,
\eqref{eq-phase-comp} shows the comparison of their derivatives at where
they take the same value,
hence we would call it the phase comparison principle.
The prototype of \eqref{eq-phase-comp} (and \eqref{eq-phase-mono})
is already formulated in the proof
of Lemma 3.6 in \cite{Non20}.
Here we omit the details.
$\hfill\Box$

The above monotone dependence is used to construct special upper and lower
solutions at speed near the critical wave speed $c^*(m,r)$.

\begin{lemma} \label{le-upper-lower}
For any $c>c^*(m,r)$, there exists a number $k>\kappa$, such that
$$
\overline u(t,x):=\overline\phi{}_c^k(\xi):=
\begin{cases}
\phi_c^k(\xi), \quad &\xi<\xi_c^k,\\
k, \quad &\xi\ge\xi_c^k,
\end{cases}
\qquad \xi=x\cdot\nu+ct,
$$
is an upper solution of \eqref{eq-main}
with the initial data $\overline u_0(s,x):=\overline\phi{}_c^k(x\cdot\nu+cs)$ for $s\in[-r,0]$,
where $\phi_c^k(\xi)$ is the unique solution of the degenerate problem \eqref{eq-tw-local}
on its maximal existence interval $(-\infty,\xi_c^k)$ with $\xi_c^k\in(0,+\infty]$.
Similarly, for any $c\in(0,c^*(m,r))$,
$$
\underline u(t,x):=\underline\phi{}_c^\kappa(\xi):=
\begin{cases}
\phi_c^\kappa(\xi), \quad &\xi<\xi_c^\kappa,\\
0, \quad &\xi\ge\xi_c^\kappa,
\end{cases}
\qquad \xi=x\cdot\nu+ct,
$$
is a lower solution of \eqref{eq-main}
with the initial data $\underline u_0(s,x):=\underline\phi{}_c^\kappa(x\cdot\nu+cs)$
for $s\in[-r,0]$,
$\xi_c^\kappa<+\infty$,
and $\sup_{\xi\in\mathbb R}\phi_c^\kappa(\xi)<\kappa$.
Moreover, for any $c\in(0,c^*(m,r))$,
$$
\hat {\underline u}(t,x):=\hat {\underline\phi}{}_c^\kappa(\xi):=
\begin{cases}
\phi_c^\kappa(\xi), \quad &\xi<\hat \xi_c^\kappa,\\
\phi_c^\kappa(\hat \xi_c^\kappa), \quad &\xi\ge\hat\xi_c^\kappa,
\end{cases}
\qquad \xi=x\cdot\nu+ct,
$$
also is a lower solution of \eqref{eq-main}
with initial data $\hat{\underline u}_0(s,x):=\hat{\underline\phi}{}_c^\kappa(x\cdot\nu+cs)$
for $s\in[-r,0]$.
Additionally, it holds
$$\lim_{c\to (c^*(m,r))^-}\sup_{\xi\in\mathbb R}\phi_c^\kappa(\xi)=\kappa.$$
\end{lemma}
{\it\bfseries Proof.}
For the critical wave speed $c=c^*$, the
unique solution of the degenerate problem \eqref{eq-tw-local} is the sharp traveling
wave $\phi^*(\xi)=\phi_{c^*}^\kappa(\xi)$,
and the phase transform function $\tilde\psi_{c^*}(\phi)$ satisfies
$\tilde\psi_{c^*}(\phi)>0$ for $\phi\in(0,\kappa)$ and
$\tilde\psi_{c^*}(0)=\tilde\psi_{c^*}(\kappa)=0$.

For any $c>c^*$, according to the phase comparison principle \eqref{eq-phase-mono}
in Lemma \ref{le-monotone}, we see that for any $k>\kappa$,
$\tilde\psi_c(\kappa)>\tilde\psi_{c^*}(\kappa)=0$,
that is, $(\phi_c^k)'((\phi_c^k)^{-1}(\kappa))>0$.
We choose $k>\kappa$ such that case (a) in Lemma \ref{le-solve} occurs, i.e.,
$\phi_c^k(\xi)$ is strictly increasing on $(0,\xi_c^k)$ and $\phi_c^k(\xi_c^k)=k>\kappa$.
Therefore, $\underline\phi{}_c^k$ is an upper solution of
the second order differential equation \eqref{eq-main}$_1$.

For any $0<c<c^*$, similar to the above analysis, according to
the phase comparison principle \eqref{eq-phase-mono}, we have
$\tilde\psi_c(\hat\phi)=0<\tilde\psi_{c^*}(\hat\phi)$ for
some $\hat\phi\in(0,\kappa)$ since $\tilde\psi_{c^*}(\kappa)=0$.
Therefore, $\phi_c^\kappa(\xi)$ is increasing up to $\hat\phi<\kappa$
and then decreases to zero, which means Case (b) in Lemma \ref{le-solve} occurs.
In this case, we have $\xi_c^\kappa<+\infty$ and $\phi_c^\kappa(\xi_c^\kappa)=0$.
Then it follows that $\underline\phi{}_c^\kappa$ is a lower solution of
the second order differential equation \eqref{eq-main}$_1$.
Furthermore, $\hat {\underline u}(t,x)=\hat {\underline\phi}{}_c^\kappa(\xi)$
is a lower solution of \eqref{eq-main} since
$\phi_c^\kappa{}'(\hat \xi_c^\kappa)=0$ at the cut-off edge.
The limit of $\lim_{c\to (c^*(m,r))^-}\sup_{\xi\in\mathbb R}\phi_c^\kappa(\xi)=\kappa$
follows from
the continuous dependence (see Lemma 3.4 in \cite{Non20} for example),
the monotone dependence Lemma \ref{le-monotone}, and the fact that
$\sup_{\xi\in\mathbb R}\phi_{c^*(m,r)}^\kappa(\xi)=\kappa$.
$\hfill\Box$

\medskip

Now we investigate the large-time evolution of the solution with semi-compact support.

\begin{lemma} \label{le-evolution-upper}
Let $u(t,x)$ be the solution of \eqref{eq-main} with the initial data
$u_0(s,x)$ semi-compactly supported and bounded
\begin{equation} \label{eq-initial-upper}
\mathrm{supp}\,u_0(s,\cdot)\subset \Pi(x_0,\nu),
\quad u_0\in L^\infty([-r,0]\times \mathbb R^n).
\end{equation}
Then
$$
\limsup_{t\to+\infty}\sup_{x\in\mathbb R^n}u(t,x)\le\kappa.
$$
\end{lemma}
{\it\bfseries Proof.}
Consider the following differential problem
\begin{equation} \label{eq-zevolution-upper}
\begin{cases}
U'(t)=-d(U)+b(U(t-r)), \quad t>0, \\
U(s)=U_0(s)\equiv \|u_0\|_{L^\infty([-r,0]\times \mathbb R^n)}, \quad s\in[-r,0].
\end{cases}
\end{equation}
The large-time asymptotic analysis of the time-delayed ordinary differential
equation \eqref{eq-zevolution-upper} shows that $\lim_{t\to+\infty}U(t)=\kappa$.
Based on the comparison principle, and taking $U(t)$ as an upper solution
of \eqref{eq-main},
we have
$$\limsup_{t\to+\infty}\sup_{x\in\mathbb R^n}u(t,x)\le
\limsup_{t\to+\infty}U(t)=\kappa.$$
The proof is completed.
$\hfill\Box$

\begin{lemma} \label{le-evolution-lower}
Let $u(t,x)$ be the solution of \eqref{eq-main} with initial data
$u_0(s,x)$ satisfying
\begin{equation} \label{eq-initial-lower}
u_0(s,x)\ge \phi_0((x-x_0)\cdot\nu),
\end{equation}
where $\phi_0(\eta)$ with $\eta=(x-x_0)\cdot\nu$
is a non-negative, continuous and non-trivial function.
Then for any compact subset $K\subset \mathbb R^n$,
$$
\liminf_{t\to+\infty}\inf_{x\in K} u(t,x)\ge\kappa.
$$
\end{lemma}
{\it\bfseries Proof.}
The proof is divided into two steps. The first one is to show that the support of
$u(t,\cdot)$ expands to including any given compact subset for large time,
and the second one is to show that the value of $u(t,x)$ within given compact subset
grows up as time increases.

Step I.
Define a Barenblatt type function
$$g(t,x)=\varepsilon(\tau+t)^{-\sigma}
\Big[\Big(\eta^2-\frac{|x-x_0|^2}{(\tau+t)^\beta}\Big)_+\Big]^d,
\quad x\in\Omega,~t\ge0,$$
where $d=1/(m-1)$, $\beta,\sigma$,
$\varepsilon$, $\eta$, and $\tau$ are positive constants, $x_0\in\mathbb R^n$.
Then by appropriately selecting $\beta$, $\varepsilon,\tau$, $\sigma$, $\eta$, and $x_0$,
the function $g(t,x)$ is
a weak lower solution of \eqref{eq-main} for all the time $t>0$. The
detailed calculations can be found in the proof of Lemma 4.4 in \cite{JDE20chemo}.
Although the value of $g(t,x)$ is decaying,
its support is expanding at a rate as $(\tau+t)^\frac{\beta}{2}$ for some $\beta>0$.
Therefore, for any given compact subsets $K_1\subset K_2\subset\mathbb R^n$,
there exists a time $t_1>0$ such that
$K_2\subset \mathrm{supp}\,u(t,\cdot)$ for any $t\ge t_1$
and $\inf_{t\in[t_1,t_2],x\in K_1}u(t_1,x)>0$ for any $t_2>t_1$.

Step II.
We assert that for any $\hat x\in K$ and any $\hat\kappa<\kappa$, there
exist a time $\hat t$ and an open neighbourhood $B(\hat x)$ of $\hat x$ such that
$u(t,x)\ge\hat\kappa$ for all $t\ge\hat t$ and $x\in B(\hat x)$.
Then the assertion $\liminf_{t\to+\infty}\inf_{x\in K} u(t,x)\ge\kappa$
follows from the finite covering theorem.
For any given $\hat\kappa<\kappa$, we define
$$\hat d(s):=d(s)+\lambda_0 s$$
with $\lambda_0>0$ sufficiently small such that
$b(s)>\hat d(s)$ for all $s\in(0,\hat \kappa]$
due to $b(s)>d(s)$ for all $s\in(0,\kappa)$.
That is, the minimal positive equilibrium for $b(s)$ and $\hat d(s)$
is located in $(\hat\kappa,\kappa)$.
Consider the following separated variable function
\begin{equation} \label{eq-zseparated}
U(t,x):=(\cos(\mu_0(x-\hat x)))_+^\frac{1}{m}\cdot g(t),
\end{equation}
with $\mu_0>0$ and function $g(t)>0$ to be determined.
We have
\begin{align*}
\Delta U^m(t,x)=-\mu_0^2(\cos(\mu_0(x-\hat x)))_+\cdot g^m(t)
=-\mu_0^2U^m(t,x),
\quad \forall x\in B_{R_0}(\hat x),
\quad \text{with~}
R_0:=\frac{\pi}{2\mu_0},
\end{align*}
and the generalized derivative (which is not Lebesgue integrable) satisfies
$$
\Delta U^m(t,x)\ge -\mu_0^2(\cos(\mu_0(x-\hat x)))_+\cdot g^m(t)
\cdot\chi_{B_{R_0}(\hat x)}(x)
\ge -\mu_0^2U^m(t,x),
$$
in the sense of distributions,
where $\chi_{B}(x)$ is the characteristic function of the set $B$.

Let us choose $g(t)\in(0,\kappa)$ and $\mu_0=\sqrt{{\lambda_0}/{\kappa^{m-1}}}$, then
$$
\Delta U^m(t,x)\ge-\mu_0^2U^m(t,x)
\ge -\mu_0^2\kappa^{m-1} U(t,x)=-\lambda_0U(t,x).
$$
In order to construct $U(t,x)$ as a lower solution of \eqref{eq-main}
for $t>t_1$ with some $T>0$, it suffices to set
\begin{equation} \label{eq-zU}
\begin{cases}
\displaystyle
\pd U t\le-\hat d(U)+b(U(t-r,x)),\quad &x\in \mathbb R^n,~t>T,
\\[1mm]
U(s,x)\le u(s,x), \quad &x\in\mathbb R^n,~s\in[T-r,T].
\end{cases}
\end{equation}
Therefore, we have
\begin{align*}
\pd U t&\le-\hat d(U)+b(U(t-r,x))\\
&=-\lambda_0U-d(U)+b(U(t-r,x))\\
&\le \Delta U^m(t,x)-d(U)+b(U(t-r,x)),
\end{align*}
which is the differential inequality in the definition of lower solutions.
As to the comparison of the initial data, we have
$$
U(s,x)\le g(s)\cdot\chi_{B_{R_0}(\hat x)}(x).
$$
According to Step I, by setting $K_1=B_{R_0}(\hat x)$
and $T=t_1+r$, $t_2=T$, we further have
$$
u(s,x)\ge \inf_{t\in[t_1,t_2],x\in K_1}u(t,x):=\varepsilon_0>0,
\quad \forall x\in B_{R_0}(\hat x), s\in[T-r,T].
$$
It follows that a sufficient condition for \eqref{eq-zU} is
\begin{equation} \label{eq-zU2}
\begin{cases}
\displaystyle
(\cos(\mu_0(x-\hat x)))_+^\frac{1}{m}\cdot g'(t)
\\
\qquad \le-\hat d((\cos(\mu_0(x-\hat x)))_+^\frac{1}{m}\cdot g(t))
+b((\cos(\mu_0(x-\hat x)))_+^\frac{1}{m}\cdot g(t-r)),\quad x\in \mathbb R^n,~t>T,
\\[1mm]
g(s)= \varepsilon_0, \quad s\in[T-r,T],
\qquad g(t)\in(0,\kappa), \quad t>T,
\end{cases}
\end{equation}
or alternatively,
\begin{equation} \label{eq-zU3}
\begin{cases}
\displaystyle
g'(t)\le \inf_{\lambda\in(0,1)}\frac{b(\lambda g(t-r))-\hat d(\lambda g(t))}{\lambda},
\quad ~t>T,
\\[1mm]
g(s)= \varepsilon_0, \quad s\in[T-r,T],
\qquad g(t)\in(0,\kappa), \quad t>T.
\end{cases}
\end{equation}
Note that
$$
\lim_{\lambda\to0^+}\frac{b(\lambda s)-\hat d(\lambda s)}{\lambda}
=b'(0)s-\hat d'(0)s=b'(0)s-(d'(0)+\lambda_0)s,
$$
which is strictly increasing for all $s>0$ since $b'(0)>d'(0)\ge0$
($\lambda_0$ is sufficiently small),
and $b(s)>\hat d(s)$ for all $s\in(0,\hat\kappa]$.
There exists a constant $\delta_0>0$ such that
$$
\inf_{\lambda\in(0,1)}\frac{b(\lambda s)-\hat d(\lambda s)}{\lambda}
\ge \delta_0>0, \quad \forall s\in[\varepsilon,\hat\kappa].
$$

We now solve the following time-delayed ordinary differential equation step by step:
\begin{equation} \label{eq-zdelay}
\begin{cases}
\displaystyle
g'(t)=\inf_{\lambda\in(0,1)}\frac{b(\lambda g(t-r))-\hat d(\lambda g(t))}{\lambda},
\quad t>T,\\
g(s)=\varepsilon_0, \quad t\in[T-r,T].
\end{cases}
\end{equation}
Firstly, for $t\in[T,T+r]$, we have
$$
g'(T)=\inf_{\lambda\in(0,1)}\frac{b(\lambda \varepsilon_0)
-\hat d(\lambda \varepsilon_0)}{\lambda}\ge \delta_0>0,
$$
which means $g(t)$ is strictly increasing
until $t\ge T+r$ or
\begin{equation} \label{eq-zlambda}
\inf_{\lambda\in(0,1)}\frac{b(\lambda \varepsilon_0)
-\hat d(\lambda g(t))}{\lambda}=0.
\end{equation}
Since there exist two constants $C_2\ge C_2>0$ such that
$C_1g(t)\le{\hat d(\lambda g(t))}/{\lambda}\le C_2g(t)$,
the asymptotic analysis of linear differential inequality shows that
\eqref{eq-zlambda} cannot happen in finite time.
Therefore, $g(t)$ is strictly increasing on $[T,T+r]$ and
\begin{equation} \label{eq-zlambda2}
\inf_{\lambda\in(0,1)}\frac{b(\lambda \varepsilon_0)
-\hat d(\lambda g(t))}{\lambda}>0, \quad \forall t\in[T,T+r].
\end{equation}
Secondly, for $t\in[T+r,T+2r]$, we have
$$
g'(T+r)=
\inf_{\lambda\in(0,1)}\frac{b(\lambda g(T))-\hat d(\lambda g(T+r))}{\lambda}
=\inf_{\lambda\in(0,1)}\frac{b(\lambda \varepsilon_0)
-\hat d(\lambda g(T+r))}{\lambda}>0,
$$
due to \eqref{eq-zlambda2}.
It follows that $g(t)$ is strictly increasing until $t\ge T+2r$ or
\begin{equation} \label{eq-zlambda3}
\inf_{\lambda\in(0,1)}\frac{b(\lambda g(t-r))
-\hat d(\lambda g(t))}{\lambda}=0,
\end{equation}
where $g(t-r)$ is already known as $t-r\in[T,T+r]$.
An asymptotic analysis shows that \eqref{eq-zlambda3} cannot happen in finite time,
especially, in $[T+r,T+2r]$.
Otherwise, let $t^*\in(T+r,T+2r]$ be the minimal time such that $\eqref{eq-zlambda3}$
is valid.
Then $g'(t)=0$ and ${b(\lambda g(t-r))}/{\lambda}$ is strictly increasing.
Hence there exists a $\hat t^*\in(T+r,t^*)$ such that
$$
\inf_{\lambda\in(0,1)}\frac{b(\lambda g(t-r))
-\hat d(\lambda g(t))}{\lambda}<0, \quad t\in(\hat t^*,t^*),
\qquad \text{and~}
\inf_{\lambda\in(0,1)}\frac{b(\lambda g(\hat t^*-r))
-\hat d(\lambda g(\hat t^*))}{\lambda}=0,
$$
which contradicts to the minimality of $t^*$.
Repeating the above arguments, we see that $g(t)$ is increasing
and the minimal positive equilibrium of
$\inf_{\lambda\in(0,1)}{(b(\lambda s)-\hat d(\lambda s))}/{\lambda}$
is greater than $\hat\kappa$.
There exists a time $\hat t>T$ such that $g(t)>\hat\kappa$ for all $t\ge\hat t$.
Furthermore, by the comparison principle,
$$
u(t,x)\ge U(t,x)=(\cos(\mu_0(x-\hat x)))_+^\frac{1}{m}\cdot g(t)
>(\cos(\mu_0(x-\hat x)))_+^\frac{1}{m}\cdot \hat\kappa, \quad t\ge\hat t.
$$
That is,
$$
u(t,\hat x)\ge U(t,\hat x)= g(t)>\hat\kappa, \quad t\ge\hat t.
$$
Based on the uniformly continuity of $(\cos(\mu_0(x-\hat x)))_+$
near $\hat x$ with respect to $t$,
we can find a neighborhood $B(\hat x)$ of $\hat x$, independent of time, such that
$u(t,x)\ge\hat\kappa$, for all $x\in B(\hat x)$ and $t\ge\hat t$.
The proof is complete.
$\hfill\Box$

Next, in order to get the large time speed of propagation, we are going to prove it by combining the large time evolution of the solution proved in Lemma \ref{le-evolution-upper},
Lemma \ref{le-evolution-lower}, and the special upper and lower solutions in
Lemma \ref{le-upper-lower}.

{\it\bfseries Proof of Theorem \ref{th-evolution}.}
First of all, from
Lemma \ref{le-evolution-upper} and Lemma \ref{le-evolution-lower}, we have the large time evolution such that
$$
\limsup_{t\to+\infty}\sup_{x\in\mathbb R^n}u(t,x)\le\kappa,
$$
and for any compact subset $K\subset \mathbb R^n$
$$
\liminf_{t\to+\infty}\inf_{x\in K} u(t,x)\ge\kappa.
$$
Note that the initial condition \eqref{eq-initial-semi} is translation invariant
in the direction perpendicular to $\nu$,
similar to the proof of Lemma \ref{le-evolution-lower},
for any finite numbers $s_1<s_2$,
there holds
$$
\liminf_{t\to+\infty}\inf_{(x-x_0)\cdot\nu\in[s_1,s_2]} u(t,x)\ge\kappa.
$$
Without loss of generality, we may assume that $\phi_0(\eta)$ is symmetric
(after shifting if necessary) with respect to $\eta=0$.
Otherwise, we can choose another function with symmetry and smaller than $\phi_0(\eta)$.
Then the lower solution $u^*(t,x)$ with the initial data given by
$u_0^*(s,x)=\phi_0((x-x_0)\cdot\nu)$
is also symmetric with respect to $\eta=0$.
The propagation properties of $u^*(t,x)$ at one side of $\eta<0$
is equivalent to the initial boundary value problem with homogeneous Neumann condition
in the half space according to the refection principle.

For any $c>c^*(m,r)$, let $\hat c\in(c^*(m,r),c)$ and
$\overline u(t,x)=\overline\phi{}_{\hat c}^k(\xi)$
with $\xi=x\cdot\nu+\hat ct$
be the upper solution of \eqref{eq-main}
corresponding to $\hat c>c^*(m,r)$ as proved in Lemma \ref{le-upper-lower}.
Note that $\lim_{\xi\to+\infty}\overline\phi{}_{\hat c}^k(\xi)=k>\kappa$,
we can change the initial time to some $T>0$ such that
$$
\sup_{x\in\mathbb R^n}u(t,x)<\frac{k+\kappa}{2}<k, \quad \forall t>T,
$$
and then shift $\overline u(t,x)=\overline\phi{}_{\hat c}^k(\xi)$
such that the comparison of the initial data is valid.
The comparison principle shows that
$$
\lim_{t\to+\infty}u(t,x-ct\nu)\le \lim_{t\to+\infty}\overline u(t,x-ct\nu)
=\lim_{t\to+\infty}\overline\phi{}_{\hat c}^k((x-ct\nu)\cdot\nu+\hat ct)
=\lim_{t\to+\infty}\overline\phi{}_{\hat c}^k(x\cdot\nu-(c-\hat c)t)=0,
$$
since $c>\hat c$ and $\overline\phi{}_{\hat c}^k(\xi)=0$ for $\xi\le\xi_0$,
where $\xi_0$ is given after the shifting.

Similarly, for any $c<c^*(m,r)$, let $\hat c\in(c,c^*(m,r))$ and
$\hat{\underline u}(t,x)=\hat{\underline\phi}{}_{\hat c}^\kappa(\xi)$
with $\xi=x\cdot\nu+\hat ct$
be the lower solution of \eqref{eq-main}
corresponding to $\hat c<c^*(m,r)$ constructed in Lemma \ref{le-upper-lower}.
Since $\sup_{\xi\in\mathbb R}\underline\phi{}_{\hat c}^\kappa(\xi)<\kappa$,
we change the initial time to some $T>0$ such that
$$
\inf_{x\cdot\nu\in[-(\hat\xi_{\hat c}^\kappa+\hat cr+1),0]} u(t,x)
>\sup_{\xi\in\mathbb R}\underline\phi{}_{\hat c}^\kappa(\xi),
\qquad t\ge T,
$$
Note that $\hat{\underline\phi}{}_{\hat c}^\kappa(\xi)=0$ for $\xi\le0$ and
$\hat{\underline\phi}{}_{\hat c}^\kappa(\xi)=\phi_{\hat c}^\kappa(\hat \xi_{\hat c}^\kappa)$
for $\xi\ge\hat\xi_{\hat c}^\kappa$,
we shift $\hat{\underline\phi}{}_{\hat c}^\kappa(\xi)$ such that
$\hat{\underline\phi}{}_{\hat c}^\kappa(\xi)=0$ for $\xi\le-(\hat\xi_{\hat c}^\kappa+1)$ and
$\hat{\underline\phi}{}_{\hat c}^\kappa{}'(\xi)=0$ for $\xi\ge-1$.
Therefore, $\hat{\underline u}(t,x)=\hat{\underline\phi}{}_{\hat c}^\kappa(\xi)$
is a lower solution of the corresponding homogeneous Neumann problem
\eqref{eq-main} on the half space.
According to the comparison principle, we have
$$
\lim_{t\to+\infty}u(t,x-ct\nu)\ge \lim_{t\to+\infty}\hat{\underline u}(t,x-ct\nu)
=\lim_{t\to+\infty}\hat{\underline\phi}{}_{\hat c}^k((x-ct\nu)\cdot\nu+\hat ct)
=\lim_{t\to+\infty}\hat{\underline\phi}{}_{\hat c}^k(x\cdot\nu+(\hat c-c)t)
=\phi_{\hat c}^\kappa(\hat \xi_{\hat c}^\kappa),
$$
since $\hat c>c$.
According to Lemma \ref{le-upper-lower},
$$\lim_{\hat c\to (c^*(m,r))^-}\sup_{\xi\in\mathbb R}\phi_{\hat c}^\kappa(\xi)
=\lim_{\hat c\to (c^*(m,r))^-}\phi_{\hat c}^\kappa(\hat \xi_{\hat c}^\kappa)
=\kappa,$$
and $\hat c\in(c,c^*(m,r))$ is arbitrary,
we see that $\lim_{t\to+\infty}u(t,x-ct\nu)\ge \kappa$.
Combining this with the fact that
$\limsup_{t\to+\infty}\sup_{x\in\mathbb R^n}u(t,x)\le\kappa$,
we have $\lim_{t\to+\infty}u(t,x-ct\nu)=\kappa$.
The proof is completed.
$\hfill\Box$

\section{Numerical simulations}

This section is devoted to the numerical simulations for the propagation properties of the
degenerate diffusion equation \eqref{eq-main} with time delay.
The most artful part is the numerical calculation of $\Delta u^m$ near the sharp edge (i.e., the boundary of the support).
For simplicity, we only consider the $1$-dimensional case.
We note that the sharp traveling wave $\phi^*(x+c^*t)$ is a typical solution that propagates to the left direction
with fixed speed $c^*>0$
and $c^*(r)<c^*(0)$ for time delay $r>0$ according to Theorem \ref{le-sharp}.

We point out that the classical second order difference scheme
\begin{equation} \label{eq-difference}
(u^m)_{xx}\big|_{x_k}\approx \frac{u^m(t,x_{k+1})+u^m(t,x_{k-1})-2u^m(t,x_k)}{(\Delta x)^2}
\end{equation}
does not work well near the boundary, since the solution of the degenerate equation has sharp moving edge.
In fact, Lemma \ref{le-solve} shows that the sharp traveling wave
\begin{equation} \label{eq-sharp-edge}
\phi^*(\xi)=\Big(\frac{(m-1)c^*}{m}\xi\Big)_+^\frac{1}{m-1}+o(|\xi|^\frac{1}{m-1}),
\quad \xi\to0.
\end{equation}
Therefore, if we take $\phi^*(x)$ as the initial data, the solution $u(t,x)$ propagates to the left direction
with the same profile.
Remarkably,  $(u^m)_{xx}$ is not continuous near the boundary for $m\ge2$,
hence the second order difference scheme based on the values at nearby discrete points loses accuracy.

{\bf Sharp-profile-based difference scheme.}
In order to handle the large variation of $(u^m)_{xx}$ near the boundary, we propose the following
sharp-profile-based difference scheme based on the expansion of the profile in \eqref{eq-sharp-edge}.
We take the case of $m=2$ for example.
Other cases can be converted to a similar equation of the so-called pressure function
$v(t,x):=\frac{m}{m-1}u^{m-1}(t,x)$.

(i) Besides the partition points $x_{-N}<x_{-N+1}<\cdots<x_0<\cdots<x_{N-1}<x_N$, we additionally introduce
an edge point $\hat x\in (x_{k-1},x_k]$ for some $-N<k\le N$ depending on time $t$ such that
$u(t,x_j)=0$ for $j<k$ and $u(t,x_j)>0$ for $j> k$;

(ii) The second order derivatives $(u^2)_{xx}\big|_{x_j}$ away from the boundary
(i.e., for $j>k$) are calculated by the classical second order difference scheme \eqref{eq-difference}
and $(u^2)_{xx}\big|_{x_j}=0$ for $j<k$ since locally $u(t,x)\equiv0$ near $x_j$;

(iii) For $(u^2)_{xx}$ near the boundary, i.e., for $(u^2)_{xx}\big|_{x_k}$, we use the profile ansatz
according to \eqref{eq-sharp-edge}
$$u(t,x)=c_1(x-\hat x)_++c_2(x-\hat x)_+^2+o((x-\hat x)_+^2),\quad x\to \hat x,$$
and the values $u(t,x_k)$, $u(t,x_{k+1})$, $u(t,x_{k+2})$ to fit the coefficients $c_1$ and $c_2$,
then
$$u^2(t,x)=c_1^2(x-\hat x)_+^2+2c_1c_2(x-\hat x)_+^3+o((x-\hat x)_+^3),\quad x\to \hat x,$$
such that we calculate
$(u^2)_{xx}\big|_{x_k}=2c_1^2+12c_1c_2(x_k-\hat x)_+$.

(iv) We compute the values $u(t+\Delta t,x_j)$ according to the equation \eqref{eq-main},
and special attention should be paid to the values $u(t+\Delta t,x_j)$ near the edge.
We use the ansatz
$$u(t+\Delta t,x)=c_1'(x-\hat x')_++c_2'(x-\hat x')_+^2+o((x-\hat x')_+^2),\quad x\to \hat x',$$
where $\hat x'$ is the new edge point,
and the new values $u(t+\Delta t,x_k)$, $u(t+\Delta t,x_{k+1})$, $u(t+\Delta t,x_{k+2})$ to fit the new coefficients $c_1'$, $c_2'$,
and the new edge $\hat x'$.
If the distance $x_k-\hat x'$ is larger than or equal to $\Delta x$, we modify
$$u(t+\Delta t,x_{k-1})=c_1'(x_{k-1}-\hat x')_++c_2'(x_{k-1}-\hat x')_+^2,$$
which means the solution propagates across $x_{k-1}$
and the values are calculated based on the edge profile
instead of the classical difference scheme based on nearby values.

\begin{figure}[htb]
\centering
\includegraphics[width=0.49\textwidth]{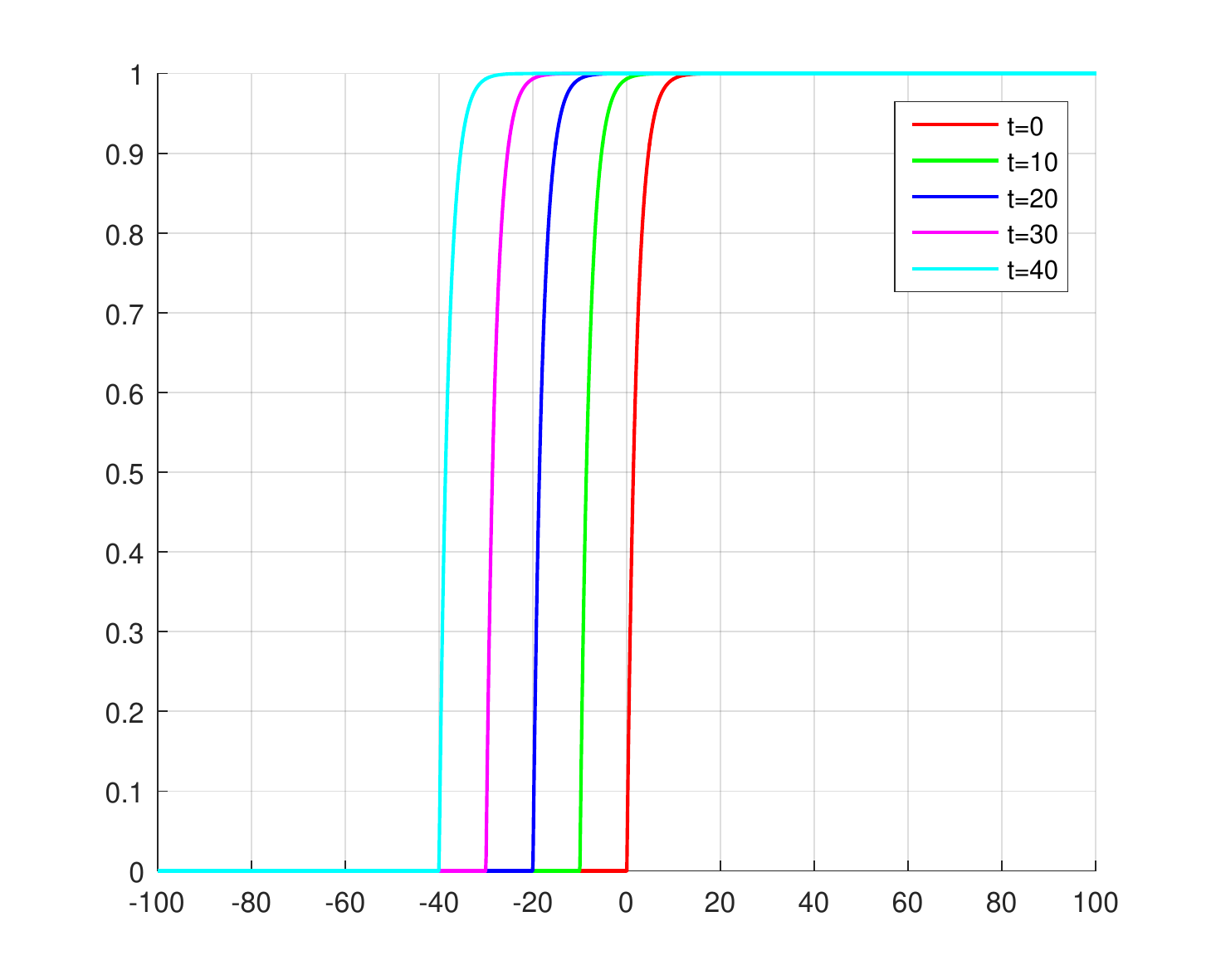}
\includegraphics[width=0.49\textwidth]{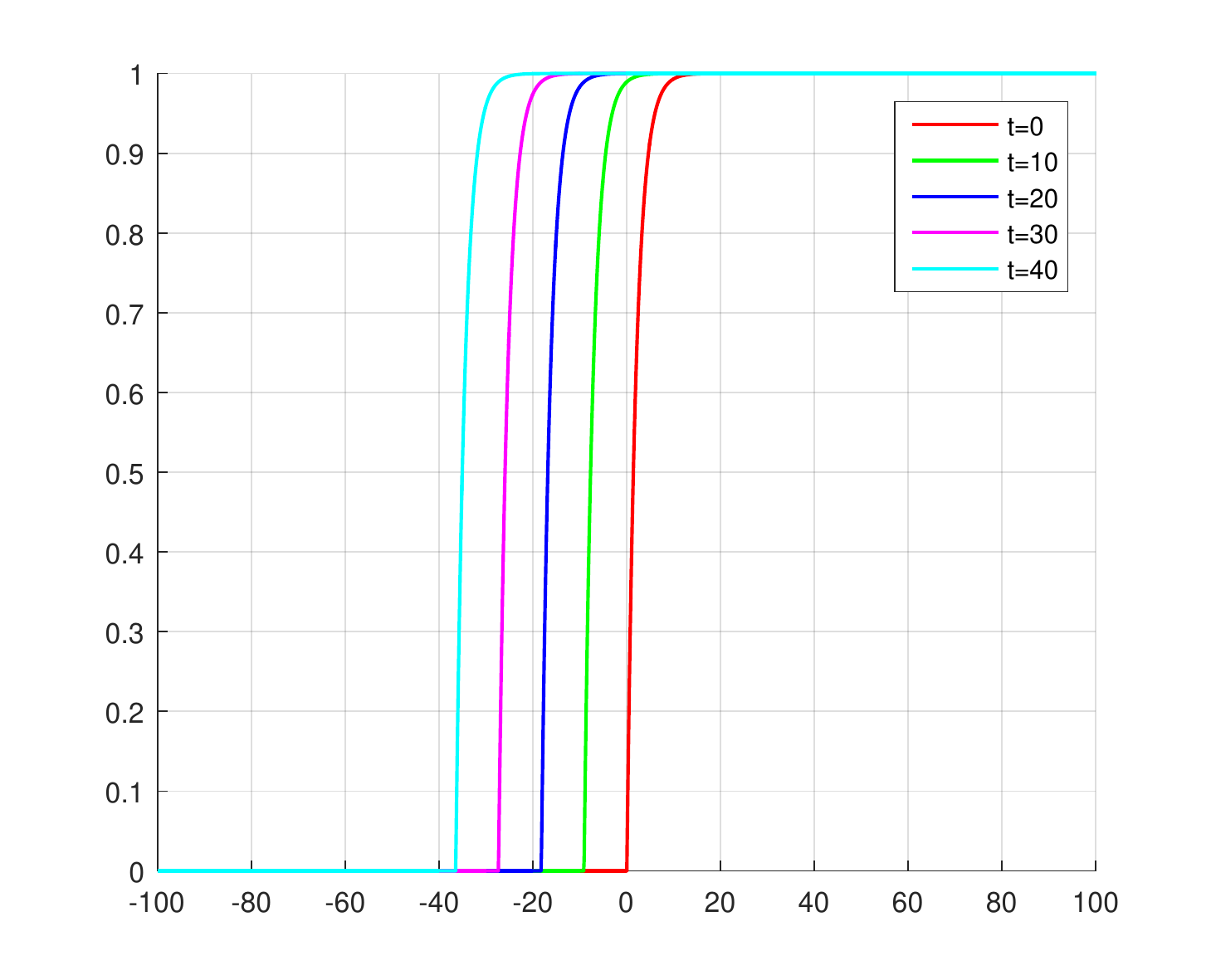}
\caption{\scriptsize The sharp-profile-based scheme for degenerate diffusion equations: (left) non-delayed case; \quad
(right) delayed case with the time delay $r=0.1$.}
\label{fig-tr}
\end{figure}

{\bf Testing the sharp-profile-based difference scheme with known profiles.}
The above sharp-profile-based scheme works perfectly for the degenerate diffusion equations with and without time delay.
For the non-delayed case, it is known that the following degenerate diffusion equation with Fisher-KPP source
\begin{equation} \label{eq-number}
u_t=(u^2)_{xx}+u-u^2,
\end{equation}
admits an explicit sharp traveling wave solution
$u(t,x)=(1-e^\frac{-x-t}{2})_+$, see \cite{Newman}.
Therefore, the minimal wave speed and the propagation speed is $c^*(2,0)=1$ for $m=2$ and $r=0$
since the sharp traveling wave is unique according to \cite{Non20}.
Note that $c^*(1,0)=2$ for the non-degenerate case $m=1$ and $r=0$.

We first take the non-smooth initial value
\begin{equation}
u_0(x):=(1-e^{-\frac{x}{2}})_+
\end{equation}
for the non-delayed degenerate diffusion equation \eqref{eq-number},
and the expected explicit solution is $u(t,x)=(1-e^\frac{-x-t}{2})_+$,
which has a sharp edge propagating to the left with speed $c^*(2,0)=1$.
The numerical simulation using the sharp-profile-based difference scheme
shows that the solution propagates with the same profile
and the sharp edge is preserved.
See the illustration Figure \ref{fig-tr}.

\begin{figure}[htb]
\centering
\includegraphics[width=0.49\textwidth]{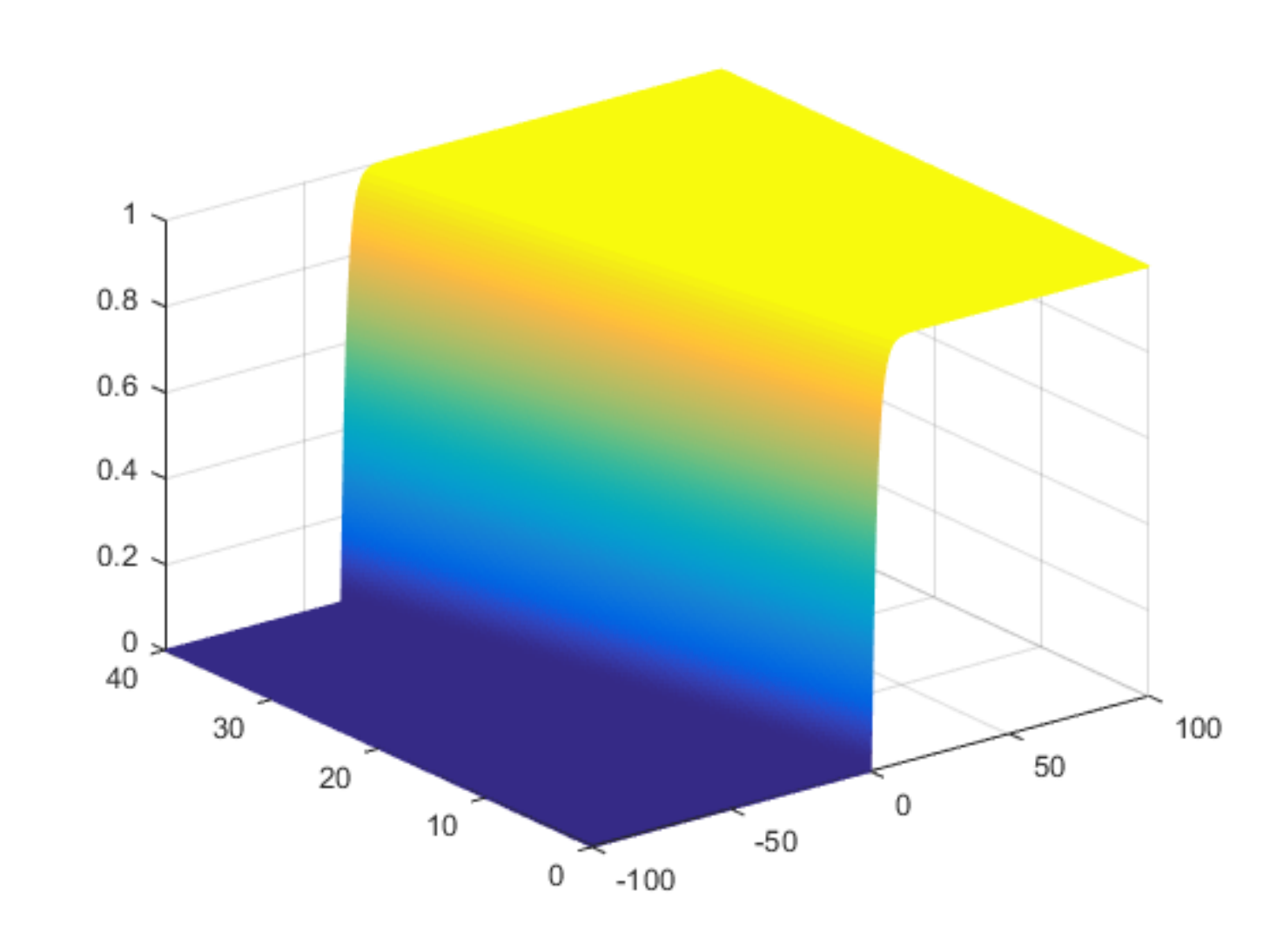}
\includegraphics[width=0.49\textwidth]{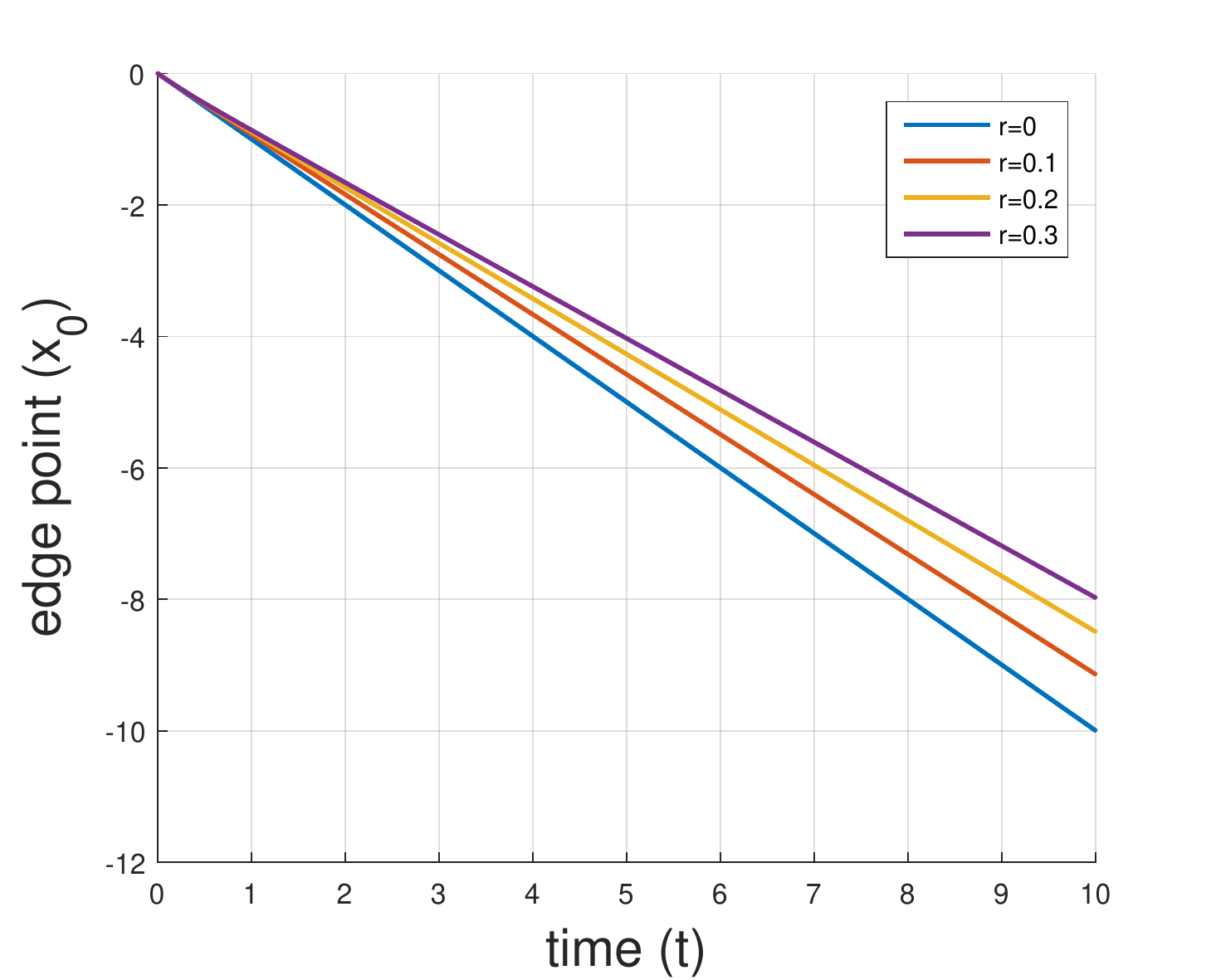}
\caption{\scriptsize (left) The $3$-dimensional image of solution for $r=0.1$; \quad
(right) The edge point as a function of time for $r=0$, $r=0.1$, $r=0.2$, and $r=0.3$, respectively.}
\label{fig-speed}
\end{figure}

{\bf Time-delayed degenerate diffusion equations.}
For the time-delayed case, we apply the above sharp-profile-based scheme to
\begin{equation} \label{eq-number2}
u_t(t,x)=(u^2(t,x))_{xx}+u(t-r,x)-u^2(t,x),
\end{equation}
with the given initial data
\begin{equation}\label{eq-number2-ID}
u_0(s,x):=(1-e^{-\frac{x}{2}})_+, \ \ \mbox{ uniformally in } s\in [-r,0].
\end{equation}
Numerical simulation shows that the solution
propagates to the left with a smaller speed $c^*(2,r)<c^*(2,0)$,
which coincides with the theoretical result in Theorem \ref{th-evolution}.
See the illustration Figure \ref{fig-tr}
and the $3$-dimensional image in Figure \ref{fig-speed}.

{\bf Propagation speed influenced by time delay.}
We take the time delay $r=0$, $0.1$, $0.2$, $0.3$, respectively.
Numerical simulations show that the corresponding propagation speed is monotonically  decreasing with respect to the time delay.
Illustrated figure is presented in Figure \ref{fig-speed}.
Here, we show that the asymptotic propagation speed after the evolution on $[0,T]$ with time $T=10$ is
\begin{equation}\label{cws-1}
c^*(2,0)\approx1.0000,  \quad \ c^*(2,0.1)\approx0.9115,  \quad \ c^*(2,0.2)\approx0.8439, \quad \  c^*(2,0.3)\approx0.7891,
\end{equation}
for the time delay $r=0, \  0.1, \  0.2$, and $ 0.3$,  respectively.

\begin{figure}[htb]
\centering
\includegraphics[width=0.49\textwidth]{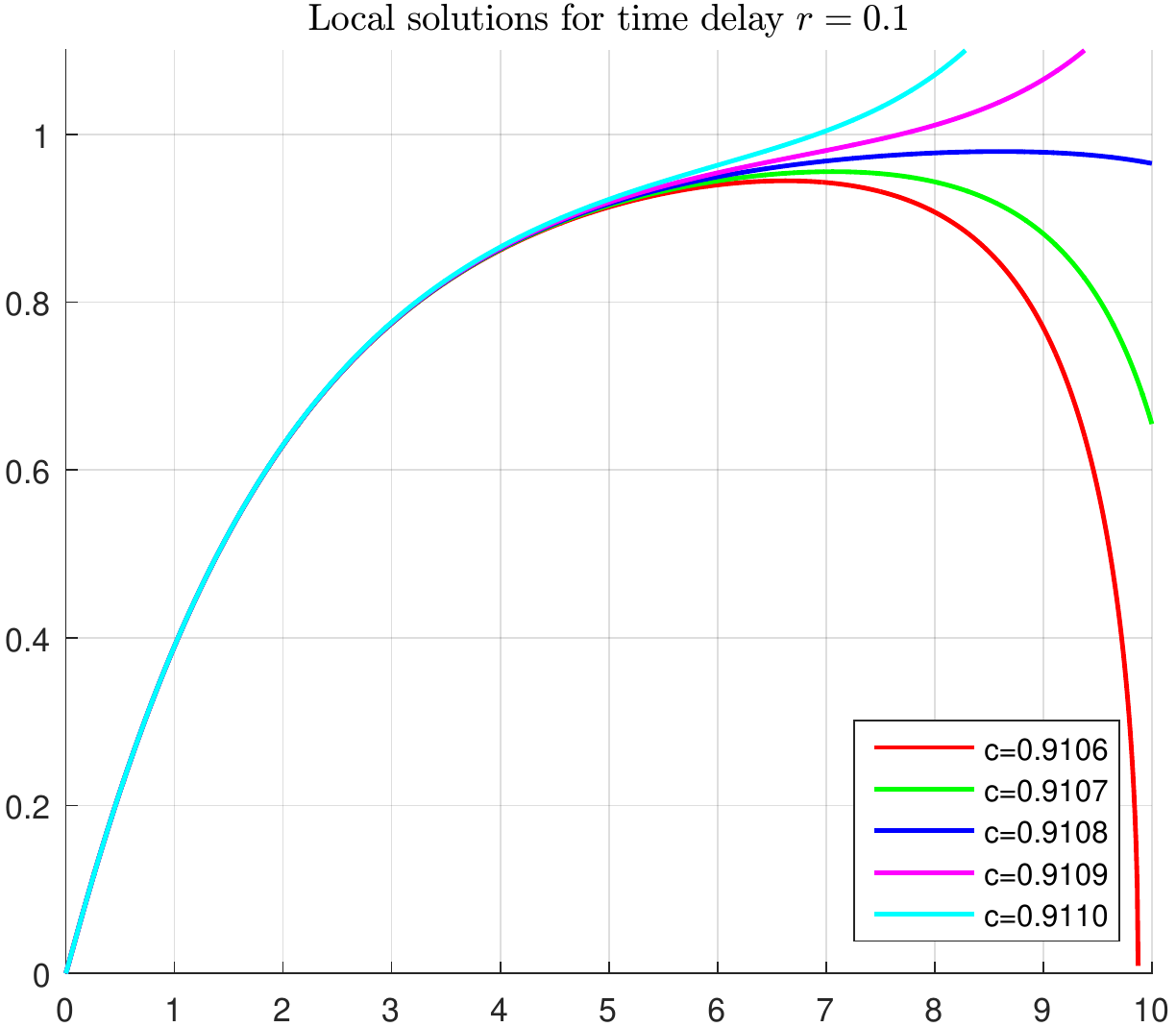}
\includegraphics[width=0.49\textwidth]{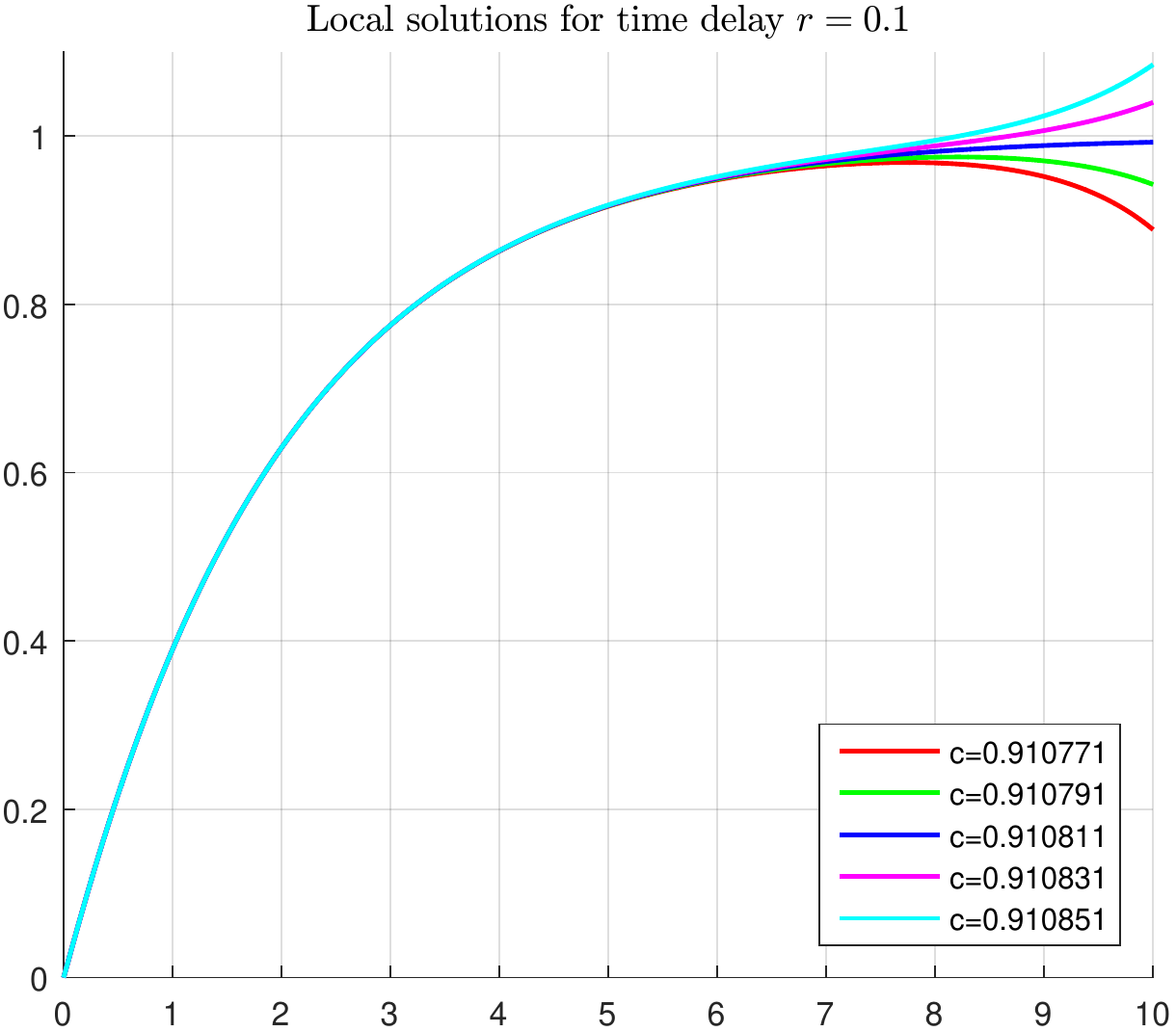}
\caption{\scriptsize The local solutions for searching of critical wave speed for $r=0.1$.}
\label{fig-cs}
\end{figure}

{\bf Coincidence with the critical wave speed for sharp traveling waves.}
The critical wave speed for the non-delayed case is know explicitly as $c^*(2,0)=1$,
while the speed of the time-delayed case is characterized by a variational inequality in \cite{Non20} such that
$c^*(m,r)<c^*(m,0)$ without information for the (numerical) calculation.
According to the proof in Section 3, the sharp wave is determined as the locally solution $\phi_c^k(\xi)$ of
a singular ODE \eqref{eq-tw-local}
such that $\phi_c^k(\xi)$ exists globally and is monotone increasing.
This primary idea of proof provides a numerical method to calculate the critical wave speed,
although it may not be efficient:\\
\indent (i) for large $c>c^*$, the local solution $\phi_c^k(\xi)$ grows up to beyond the positive equilibrium $1$;\\
\indent (ii) for small $c<c^*$, the local solution $\phi_c^k(\xi)$ eventually declines down to $0$;\\
\indent (iii) the local solution $\phi_c^k(\xi)$ is monotonically increasing with respect to $c$.\\
Based on the above theoretical observations, we carry out the following simulation of the critical wave speed
for the time delay $r=0.1$ in Figure \ref{fig-cs}.
This kind of simulation shows that
\begin{equation}\label{cws-2}
c^*(2,0)\approx 1.0000, \quad \ c^*(2,0.1)\approx0.9108, \quad \ c^*(2,0.2)\approx0.8430, \quad  \ c^*(2,0.3)\approx0.7880,
\end{equation}
for $r=0, \ 0.1$, $0.2$, and $0.3$, respectively.
We see that the numerical propagation speeds shown in \eqref{cws-1} coincide with these calculated
 critical wave speeds \eqref{cws-2}, and the numerical errors are somewhat like $10^{-4}$, which are pretty good.

\begin{figure}[htb]
\centering
\includegraphics[width=0.49\textwidth]{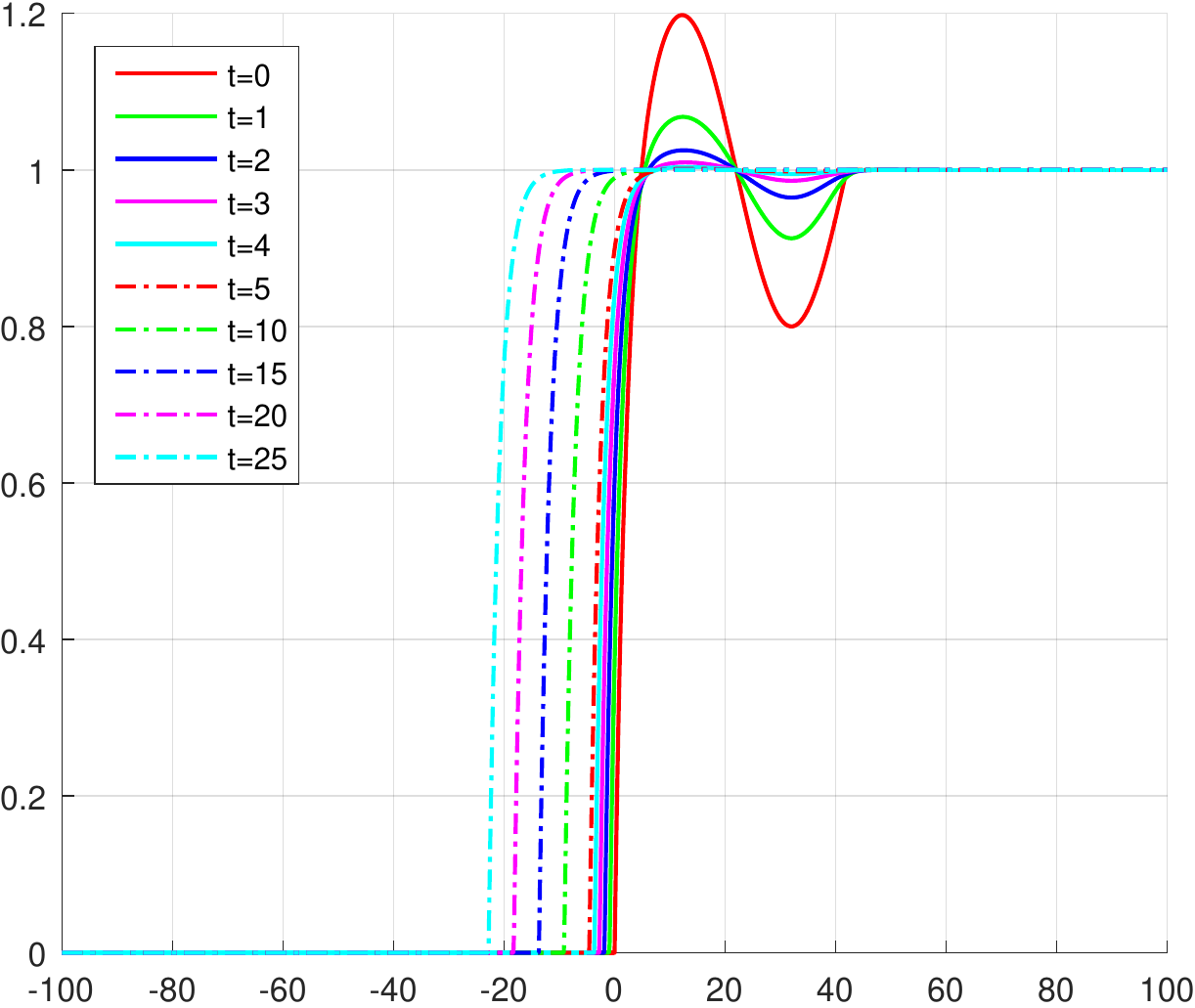}
\includegraphics[width=0.49\textwidth]{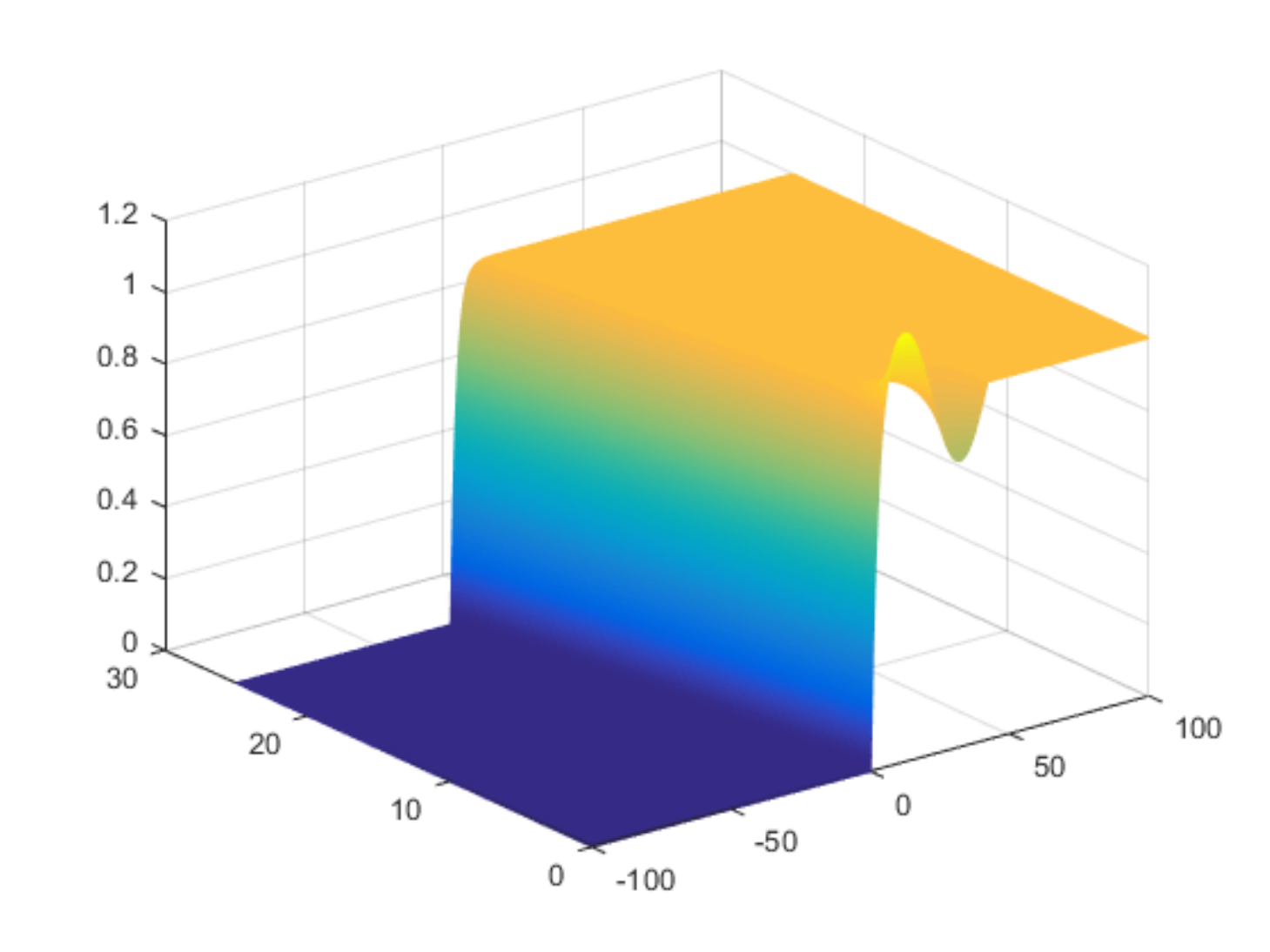}
\caption{\scriptsize The evolution of perturbation: (left) the $2$-dimensional image; \quad
(right) the $3$-dimensional image.}
\label{fig-pert}
\end{figure}

{\bf Numerical stability.}
We take the initial value
\begin{equation}\label{eq-number-2-ID-new}
u_0(s,x):=(1-e^{-\frac{x}{2}})_++0.2\sin(\pi(x-2)/20)\cdot\chi_{[2,42]}(x), \ \ \mbox{ uniformlly in } s\in [-r,0],
\end{equation}
where the perturbation is chosen within the support of the profile.
According to the numerical simulation in Figure \ref{fig-pert},
we see that the perturbation decays to zero, and the solution still propagates to the left with sharp edge.

\medskip

{\bf Acknowledgement}.
The research of S. Ji was supported in part by
by NSFC Grant No. 11701184,
Guangdong Basic and Applied Basic Research Foundation,
and the Fundamental Research
Funds for the Central Universities of SCUT.
The research of M. Mei was supported in part by
NSERC Grant RGPIN 354724-16, and FRQNT Grant No. 2019-CO-256440.
The research of T. Xu was supported in part by NSFC Grant No. 11871230 and NSFC Grant No. 11771156.
The research of J. Yin was supported in part by NSFC Grant No. 11771156,
NSF of Guangzhou Grant No. 201804010391,
and Guangdong Basic and Applied Basic Research Foundation Grant No.2020B1515310013.

\end{document}